\documentclass[graybox]{svmult}

\usepackage{type1cm}        
\usepackage{makeidx}         
\usepackage{graphicx}        
\usepackage{multicol}        
\usepackage[bottom]{footmisc}

\usepackage{newtxtext}       %

\makeindex             

\usepackage{enumerate,amsmath,amssymb,stmaryrd} 
\usepackage{pstricks,pst-node,pst-coil,pst-plot} 
\usepackage{hyperref}
\usepackage{esint}
\usepackage{xfrac}
\usepackage{accents}

\newcommand{\definedas}{\mathrel{\raise.095ex\hbox{\rm :}\mkern-5.2mu=}}
 \newcommand{\asdefined}{\mathrel{=\mkern-5.2mu\raise.095ex\hbox{\rm :}}}

\def\epsilon{{\varepsilon}}
\def\eps{{\varepsilon}}
\def\phi{{\varphi}}
\newcommand{\IR}{\mathbb{R}}
\newcommand{\qtext}[1]{\quad\text{#1}\quad}
\newcommand{\stext}[1]{\ \text{#1}\ }

\let\<\langle 
\let\>\rangle

\newcommand{\R}{\mathbb{R}}

\renewcommand{\hbar}{\overline{h}}

\DeclareMathOperator{\tr}{tr}

\newcommand{\scal}{\mathrm{Scal}}

\newcommand{\diver}{\operatorname{div}}

\newcommand{\odd}{\textrm{odd}}
\newcommand{\even}{\textrm{even}}

\newcommand{\SigST}{\Sigma_\sigma^{\text{STCMC}}}
\newenvironment{enum}{\begin{enumerate}[a)]}{\end{enumerate}}

\begin{document}
\title*{Geometrically defined asymptotic coordinates\\ in General Relativity}
\author{Carla Cederbaum\orcidID{0000-0002-3969-4850} and\\ Jan Metzger\orcidID{0000-0001-5632-7916}}
\institute{Carla Cederbaum \at T\"ubingen University, Fachbereich Mathematik, Auf der Morgenstelle 10, 72076 T\"ubingen, Germany\\ \email{cederbaum@math.uni-tuebingen.de}
\and Jan Metzger \at Potsdam University, Institut f\"ur Mathematik, Karl-Liebknecht-Str. 24-25,
14476 Potsdam-Golm, Germany\\ \email{jan.metzger@uni-potsdam.de}}

\maketitle

\abstract{We review and announce recent results on the
  asymptotic behavior of asymptotically Euclidean relativistic initial
  data sets and asymptotic foliations thereof. In particular, we discuss the
  geometrization of asymptotic flatness and of asymptotic geometric
  (in-)\-variants such as mass, energy, linear momentum, angular
  momentum and center of mass as well as their relations to certain
  geometric asymptotic foliations such as the CMC- and
  STCMC-foliations.}
  
  \keywords{Relativistic initial data, asymptotic flatness, foliations, constant mean curvature, mass, energy, momentum, center of mass}\\

\noindent\textbf{MSC codes:} 35-02, 53-02, 58-02, 35J05, 35Q75, 53A10, 53C12, 53C80
  
\section{Introduction}\label{sec:intro}
In general relativity, one central object of study are \emph{initial data sets $(M,g,K)$} consisting of a smooth Riemannian $3$-manifold $(M,g)$ and a symmetric tensor field $K\in\Gamma(T^{0}_{2}M)$ to be thought of as the extrinsic curvature (second fundamental form) of the spacelike hypersurface $M$ with induced metric $g$ sitting inside some ambient \emph{spacetime} (aka time-orientable Lorentzian manifold). If the ambient spacetime models a relativistic system isolated from external influences, it is useful and customary to restrict one's attention to its so-called ``asymptotically Euclidean'' initial data sets, modeling spacelike hypersurfaces in the spacetime obtained as the $\{t=\text{const}\}$ level sets of some asymptotically freely falling observer with time function $t$. Typically, an initial data set $(M,g,K)$ is called \emph{asymptotically Euclidean} (or \emph{asymptotically flat} or \emph{AE} for short) if, outside some compact set $C\subset M$, there is a smooth diffeomorphism
\begin{align}
\Phi\colon M\setminus C\to\R^{3}\setminus\overline{B}
\end{align}
onto the complement of a closed ball $\overline{B}\subset\R^{3}$ and if, in the Cartesian coordinates $(x^{i})$ induced by $\Phi$, one has
\begin{align}\label{eq:gdecay}
(\Phi_{\ast} g)_{ij}&=\delta_{ij}+{O}_{2}(r^{-\tau}),\\\label{eq:Kdecay}
(\Phi_{\ast} K)_{ij}&=\phantom{\delta_{ij}+}\;{O}_{1}(r^{-(\tau+1)})
\end{align}
for $i,j=1,2,3$ as $r\to\infty$ for some \emph{decay rate
  $\tau>\frac{1}{2}$}, where $r$ denotes the standard radial
coordinate on $\R^{3}$, $\delta$ denotes the Euclidean metric, and
$\mathcal{O}_{k}(r^{\alpha})$ denotes the class of smooth functions in
the $r^{\alpha}$-weighted function space
$C^{k}_{\alpha}(\R^{n}\setminus\overline{B})$ (or in the corresponding
weighted Sobolev space, see e.g.~\cite{Bartnik86}).

For an initial data set $(M,g,K)$, the Hamiltonian constraint scalar $\rho$ (aka the \emph{energy density}) and the
momentum constraint one-form $J$ (aka the \emph{momentum density}) are defined by  
\begin{equation*}
  16\pi\rho \definedas \scal_{g}-\vert K\vert_{g}^{2}+(\tr_{g}K)^{2}
  \qquad\text{and}\qquad  
    8\pi J  \definedas \operatorname{div}_{g}\left(K-(\tr_{g}K)g\right).
\end{equation*}
Here $\scal_{g}$ denotes the scalar curvature of $g$ and
$\vert\cdot\vert_{g}$, $\tr_{g}$, and $\operatorname{div}_{g}$ denote the
tensor norm, trace, and divergence with respect to $g$, respectively.

In addition to the above decay assumptions on $g$ and $K$, one also
requires the physically motivated integrability assumptions
\begin{align}\label{eq:Hamint}
  \Phi_{\ast} \rho &\in L^{1}(\R^{3}\setminus\overline{B}), \\
  \label{eq:momint}
  (\Phi_{\ast}J)_{i} &\in L^{1}(\R^{3}\setminus\overline{B})
  \qtext{for} i=1,2,3.
\end{align}
The coordinates $(x^{i})$ induced by any such diffeomorphism $\Phi$
are called \emph{asymptotic(ally Euclidean) coordinates} for the
initial data set $(M,g,K)$. For the remainder of this article, we will
suppress the push forwards in the above definitions and write
e.g. $g_{ij}$ in place of $(\Phi_*g)_{ij}$.

It might appear counter-intuitive to geometers to define asymptotic flatness of initial data sets via the existence of suitable asymptotic coordinates; instead, one might wish to \textbf{characterize suitable asymptotic behavior by more geometric (i.e., coordinate invariant) concepts}. This turns out to be anything but straightforward and part of this article is dedicated to the discussion of promising approaches to such a geometric definition, see Section~\ref{sec:geometric}.

The threshold $\tau>\frac{1}{2}$ on the decay rate $\tau$ in the decay
assumptions \eqref{eq:gdecay}, \eqref{eq:Kdecay} is chosen as small as
possible to still allow to define the physical quantities ``(total)
mass, energy, and linear momentum'' as geometric invariants of
asymptotically Euclidean initial data sets, see
Section~\ref{sec:InitialData}. 

In a similar spirit, it is also desirable to study the physical quantities \textbf{(total) angular momentum and center of mass} of asymptotically Euclidean initial data sets. This turns out to be significantly more subtle as one cannot expect these quantities to be adequately represented by geometric invariants on the basis of their physical properties such as their transformation behavior. In the approaches pursued in the literature, this non-invariance leads to convergence issues which are usually handled by assuming additional asymptotic properties such as the asymptotic parity conditions by Regge--Teitelboim or a refinement of \eqref{eq:gdecay}, \eqref{eq:Kdecay} for $\tau>1$, replacing $\delta_{ij}$ by the components of the so-called ``spatial Schwarzschild metric'' . In Sections \ref{sec:CoM} and \ref{sec:foliation}, we will discuss \textbf{why none of these approaches is fully satisfactory and which ideas might be more promising}. This will lead us directly to the question of geometrically characterizing asymptotic flatness to be discussed in Section \ref{sec:geometric}.

Related results on the existence of local and asymptotic foliations with certain geometric properties will be discussed along the way.



\section{Physical quantities as geometric invariants of initial data sets}
\label{sec:InitialData}

In 1959, Arnowitt, Deser, and Misner~\cite{ADM:1959} \emph{(ADM)} defined the \emph{(total) energy $E$ and linear momentum $\vec{P}=(P^{1},P^{2},P^{3})^{t}$} of asymptotically Euclidean initial data sets $(M,g,K)$ (with decay rate $\tau=1$) via the surface integral formulas
\begin{align}\label{def:E}
E&:=\frac{1}{16\pi}\lim_{r\to\infty}\sum_{i,j=1}^{3}\int_{\mathbb{S}^{2}_{r}}\left(g_{ij,i}-g_{ii,j}\right)\tfrac{x^{j}}{r}dA_{\delta},\\\label{def:P}
P^{j}&:=\frac{1}{8\pi}\lim_{r\to\infty}\sum_{i=1}^{3}\int_{\mathbb{S}^{2}_{r}}\left(K_{ij}-(tr_{g}K)g_{ij}\right)\tfrac{x^{i}}{r}dA_{\delta}.
\end{align}
Recall that we are suppressing the push-forward by $\Phi$, ``$_{,i}$'' denotes an $i$-th partial derivative, $\frac{\vec{x}}{r}$ denotes the Euclidean unit normal to  and $dA_{\delta}$ denotes the Euclidean area element induced on the coordinate sphere $\mathbb{S}^{2}_{r}$. The \emph{(total) mass $M$} is then defined as the Minkowskian length of the \emph{energy-momentum $4$-vector $(E,\vec{P})^{t}$}, that is, by 
\begin{align}\label{def:M}
M&:=\sqrt{E^{2}-\vert\vec{P}\vert^{2}},
\end{align}
with $\vert\cdot\vert$ denoting the Euclidean norm --- provided the
radicand is non-negative. In 1983, Denisov and Soloviev explored the
effect of the decay rate $\tau$ in \eqref{eq:gdecay} on $E$ in \eqref{def:E} and found, expressed in modern terms, that
$\tau\leq\frac{1}{2}$ allows for asymptotic coordinates on Euclidean
space in which $E\neq0$, a physically undesirable effect. This
motivated the mathematical study of establishing $E$ as a geometric
invariant (i.e., well-defined and independent of the choice of
asymptotic coordinates) when $\tau>\frac{1}{2}$ achieved by
Bartnik~\cite{Bartnik86} using asymptotic coordinates which are
harmonic with respect to $g$. Bartnik's proof significantly uses the integrability condition
\eqref{eq:Hamint}. 

Similar ideas relying on the application of the
divergence theorem to the integrability condition \eqref{eq:momint}
assert that the vector $\vec{P}$ is a geometric invariant in the
following sense: First, its components are well-defined in any
asymptotic coordinate system (with $\tau>\frac{1}{2}$). Second, if $\tau\not\in\mathbb{N}$, it can be shown \cite{Bartnik86}
that two asymptotic coordinate systems $x=(x^{i})$,
$y=(y^{i})$ are related by
\begin{align}\label{eq:transfo}
{y}&=Q{x}+O_{3}(r^{1-\tau})
\end{align}
as $r=\vert{x}\vert\to\infty$ for some orthogonal transformation $Q$. Then the corresponding \emph{momentum vectors $\vec{P}^{x}$ and  $\vec{P}^{y}$} are related by
\begin{align}\label{eq:transfoP}
\vec{P}^{{y}}&=Q^{t}\vec{P}^{{x}},
\end{align}
see Chru\'sciel~\cite{Chrusciel}.
\begin{remark}[Asymptotic Poincar\'e equivariance]
More generally, it was shown by Chru\'sciel~\cite{Chrusciel} that the energy-momentum $4$-vector $(E,\vec{P})^{t}$ is a geometric invariant of the ambient spacetime for $\tau>\frac{1}{2}$ in a similar sense, namely with Lorentz transformations replacing rotations. This result is apparently established with respect to one fixed asymptotic coordinate system, a technicality which can be removed, see forthcoming work by Cederbaum~\cite{Csuper}.
\end{remark}

\begin{remark}[Spatial super-translations]\label{rem:supertrans}
We would like to note that the non-leading order term in \eqref{eq:transfo} can be growing when $\tau<1$ and in fact can grow logarithmically even when $\tau=1$. This motivates the terminology that the coordinates $(y^{i})$ can be \emph{(spatially) super-translated} relative to the coordinates $(x^{i})$ upon recalling that a coordinate translation would take the form ${y}={x}+{a}$ for some ${a}\in\R^{3}$. 
 \end{remark}
 
 \begin{remark}[Higher dimensions]
Similar considerations can be made in $n>3$ dimensions and one obtains the threshold $\tau>\frac{n-2}{2}$ from the same reasoning.
\end{remark}

Defining asymptotic flatness via asymptotic decay, we can now of
course ask ourselves which decay rates are physical or otherwise
natural to consider. Clearly, by the above, it seems to be advisable
to assume $\tau>\frac{1}{2}$. On the other hand, it is easy to compute
that $E=\vec{P}=0$ when $\tau>1$ which suggests for physical reasons
that one should focus on $\frac{1}{2}<\tau\leq1$. It is now natural to
ask whether decay of rates $\frac{1}{2}<\tau<1$ actually occurs or
whether it is due to a ``bad'' choice of coordinates. By choosing
appropriate spacelike hypersurfaces, that is,s a time coordinate, in the
Minkowski (or Schwarzschild, or Kerr) spacetime, all these decay rates
can actually be produced with respect to the standard spatial
coordinates. These low decay rates occur in $K$, not in $g-\delta$,
see~\cite{CederbaumGrafMetzger2023}. Moreover, for any
$\frac{1}{2}<\tau<1$, Cederbaum and Maxwell~\cite{CM} are constructing 
initial data sets $(\R^{3},g,K=0)$ which are asymptotically Euclidean
of prescribed order $\frac{1}{2}<\tau<1$ but no better, even upon changing asymptotic
coordinates as in \eqref{eq:transfo}, see Section \ref{sec:AS} for related results.


\section{Center of mass and Regge--Teitelboim conditions}
\label{sec:CoM}
The ADM-definitions of energy~\eqref{def:E} and linear momentum~\eqref{def:P} were derived from a Hamiltonian consideration of the Einstein equations. This derivation was
later extended by Regge and Teitelboim
\cite{Regge--Teitelboim-1974} in order to define (total) angular momentum --- see Section~\ref{sec:ang} --- and center of mass. Based on their analysis, Beig and \'O
Murchadha \cite{Beig-OMurchadha-1987} define the \emph{Beig--\'O Murchadha--Regge--Teitelboim center of mass} or
\emph{B\'ORT-center of mass} $\vec{Z}_\text{B\'ORT}$ with components
\begin{equation}
  \label{eq:BORTcenter}
  Z^i
  \definedas \frac{1}{16\pi E}\lim_{r\to\infty} \sum_{k,l=1}^{3}\! \int_{\mathbb{S}^{2}_r} \!\left[ x^i \big(g_{kl,k} - 
    g_{kk,l}\big) \tfrac{x^l}{r} - \big(g_{ik} \tfrac{x_k}{r} - g_{kk}
    \tfrac{x^i}{r}\big)\right] \! dA_\delta,
\end{equation}
$i=1,2,3$, provided that $E\neq0$. The same expression is obtained by Michel~\cite{Michel} using a geometric approach via Killing initial data. Note that this limit need not exist if one just assumes asymptotic flatness, i.e., \eqref{eq:gdecay}--\eqref{eq:momint}, since the integrand generically only decays as $O(r^{-\tau})$ as $r\to\infty$.

This divergence can be remedied by imposing asymptotic parity conditions
suggested by Regge and Teitelboim \cite{Regge--Teitelboim-1974}. To
state these, we fix asymptotic coordinates $(x^{i})$ on a given initial data set $(M,g,K)$ and, continuing to suppress the diffeomorphism $\Phi$, we denote the \emph{odd part of $g$} by
\begin{equation*}
  g^\odd_{ij}(x) \definedas  \tfrac12\big(g_{ij}(x) - g_{ij}(-x)\big)
  \quad\text{for}\quad
  x\in \R^3 \setminus \overline{B}.
\end{equation*}
Similarly, the \emph{even part of $K$} is defined by 
\begin{equation*}
  K^\even_{ij}(x) \definedas  \tfrac12\big(K_{ij}(x) + K_{ij}(-x)\big)
  \quad\text{for}\quad
  x\in \R^3 \setminus \overline{B}.
\end{equation*}
Furthermore we need
\begin{equation*}
  \rho^\odd (x) \definedas  \tfrac12 \big(\rho(x) - \rho(-x)\big)
  \qquad\text{and}\qquad
  J_i^\odd (x) \definedas  \tfrac12 \big(J_i(x) - J_i(-x)\big).
\end{equation*}
Using this notation, we then say that the initial data set $(M,g,K)$ \emph{satisfies the Regge-Teiltelboim
conditions with decay rates
$(\tau,\eta) \in (\tfrac12,1] \times (0,1]$} if $(M,g,K)$ is
asymptotically flat with decay rate $\tau$ as in
Section~\ref{sec:InitialData} and if, in addition,
\begin{equation*}
 g^\odd_{ij} = O_2(r^{-1-\tau-\eta})
  \qtext{and}
  K^\even_{ij} = O_1(r^{-2-\tau-\eta})
\end{equation*}
hold as $r\to\infty$ and the integrability conditions
\begin{equation}
  \label{eq:RTint}
   x \rho^\odd \in L^1(\R^3\setminus \overline{B})
  \qtext{and}
   x J_i^\odd \in L^1(\R^3 \setminus\overline{B})
\end{equation}
are satisfied componentwise for $i=1,2,3$. We refer to the case $\eta=\tfrac12$ as the \emph{weak Regge--Teitelboim
conditions} and to the case $\eta=1$ as the \emph{strong Regge--Teitelboim
conditions}. This is because $\eta=\frac{1}{2}$ and $\eta=1$ represent thresholds in the analysis, while the remaining decay can be incorporated into the decay rate $\tau$.

The significance of the Regge--Teitelboim conditions stems from the
fact that the strong Regge--Teitelboim conditions guarantee the
convergence of the B\'ORT-center of mass integral in
\eqref{eq:BORTcenter} and similar convergence properties of the Regge--Teitelboim angular momentum, see Section~\ref{sec:ang}. To get an idea of why this holds true, note
that the coordinate spheres $\mathbb{S}^{2}_{r}$ are point reflection
symmetric across the origin so that each surface integral in~\eqref{eq:BORTcenter} only
depends on the components of $g^\text{odd}$ (modulo terms that will not matter for the
limit). Arguing in a similar spirit as Bartnik~\cite{Bartnik86} in his
proof of $E$ being a geometric invariant, the surface integrals
in~\eqref{eq:BORTcenter} can then be written as a volume integral of a
divergence which equals $x^i \rho$ up to terms which decay fast enough
to be integrable. Assumption~\eqref{eq:RTint} then guarantees the
integrability of the divergence term, an idea first explored in \cite{CorvinoWu}.

 Note that an important feature here is that the point reflection symmetry of the whole setup allows to 
require stronger decay only on $g^\odd$ and $K^\even$, which allows
initial data $(M,g,K)$ to satisfy the strong Regge--Teitelboim conditions
\emph{and} to have non-zero energy and linear momentum --- in contrast to the naive approach of raising $\tau>1$ discussed in Section~\ref{sec:InitialData}.

Moreover, it is important to note that the definition of odd and even parts of functions (metric components etc.) given above subtly depends on a choice of asymptotic coordinates. This leads to subtly strange transformation properties of the Regge--Teitelboim conditions under coordinate changes, see our forthcoming work with Graf \cite{CederbaumGrafMetzger2023} and Section~\ref{sec:noRT}. 

\subsection{Are the Regge--Teitelboim conditions generic?}\label{sec:noRT}
It should maybe not be surprising that the asymptotic point reflection symmetry imposed by the (weak or
strong) Regge--Teitelboim conditions is actually restrictive: One of the results of this project is that there do exist initial data sets $(M,g,K)$ which do \emph{not} carry any such coordinates, see \cite{Cederbaum2023} as well as the forthcoming work by Cederbaum, Graf, and Metzger~\cite{CederbaumGrafMetzger2023}:

First of all, it is not too difficult to construct initial data sets
$(M,g,K)$ which do not satisfy the Regge--Teitelboim conditions in a
particular asymptotic coordinate system which is provided from the
construction itself. To be specific, following ideas of Cederbaum and
Nerz \cite{Cederbaum-Nerz-2015}, one can consider graphical spatial slices of
the Schwarzschild spacetime: Denote
by $I_m \definedas  (\max\{0,2m\},\infty) \subset (0,\infty)$. Then the
Schwarzschild spacetime of mass $m\neq0$ is given by
$(\overline{M}_m \definedas  \IR \times M_m, \overline{g}_m)$, where
$M_m\definedas  I_m \times \mathbb{S}^2\subseteq\R^3$. We denote the coordinate along the
$\IR$-factor in $\overline{M}_m$ by $t$ and the coordinate along the
interval $I_m$ by $r$. The metric $\overline{g}_m$ is then given by
\begin{equation}
  \label{eq:Schwarzschild-spacetime}
  \overline{g}_m = -N_{m}^2 dt^2 + g_m, \stext{with}
  g_m = N_m^{-2} dr^2 + r^2 d\Omega^2, 
  N_m(r) = \sqrt{1-\frac{2m}{r}}.
\end{equation}
The desired examples arise from graphical slices of $\overline{M}_m$ of the form
\begin{equation*}
M_f \definedas  \{t=f(x) \mid x\in M_m\},
\end{equation*}
where we consider graph functions
$f = f_\beta^\even + f_\gamma^\odd$ with
\begin{equation*}
  f_\beta^\even(r) =\begin{cases} \sin(\ln(r)) & \stext{if} \beta =0 \\
    r^{\beta} & \stext{if} \beta\in(-\infty,\tfrac12)
  \end{cases}
\end{equation*}
and 
\begin{equation*}
f^\odd(x) = \frac{\< x, a \>}{r^{1-\gamma}},
\end{equation*}
 where $x$ denote the Cartesian coordinates constructed from the radial coordinate $r$ and polar coordinates on $\mathbb{S}^{2}$, $a\in\R^{3}$ is fixed with $a\neq0$, and $\<\cdot,\cdot\>$ denotes the Euclidean inner product. The case $\gamma=0$
and $\beta=0$ reduces to the examples of Cederbaum and Nerz~\cite{Cederbaum-Nerz-2015} discussed in
Section~\ref{sec:foliation}. Different choices of $\beta$ and $\gamma$
yield a range of examples which do not satisfy the weak or strong
Regge--Teitelboim conditions in the $x$-coordinates induced on $M_f$.

Our main contribution in \cite{CederbaumGrafMetzger2023} is to show
that these explicit initial data do not satisfy the weak or strong Regge--Teitelboim
conditions in \emph{any} asymptotically flat coordinate system. Our main tool here is to analyze the transformation of the
$x$-coordinates to the nearby harmonic asymptotic coordinates constructed by
Bartnik~\cite{Bartnik86}. In the harmonic coordinates $x_h$
constructed from $x$, the initial data sets $(M_{f},g_{f},K_{f})$ also do not satisfy
the Regge--Teitelboim conditions.

On the other hand, if $(M_f,g_{f},K_{f})$ would have an asymptotically
Euclidean coordinate system $y$ in which the Regge--Teitelboim
conditions hold, then the harmonic coordinates $y_h$ constructed
from $y$ also have the property that $(M_{f},g_{f},K_{f})$ satisfies the
Regge--Teitelboim conditions in these coordinates\footnote{This statement is actually not very precise, since we have to a) do the analysis in a suitable H\"older space and this require slightly more than the Regge--Teitelboim conditions stated above and b) we lose a derivative and an arbitrary small $\eps>0$ of decay rate in the Regge--Teitelboim conditions when changing coordinates from $x$ to $x_h$.}.  Since Bartnik showed that harmonic coordinates are
unique up to asymptotic rotations and translations, in fact $y_h$ is
related to $x_h$ by such an asymptotic rotation and translation. Such
transformations essentially preserve the Regge--Teitelboim conditions and we arrive
at the desired contradiction.

\subsection{An alternative to the Regge--Teitelboim conditions}\label{sec:AS}
An alternative to explicitly assuming Regge--Teitelboim conditions is to ask for asymptotic spherical symmetry, namely that 
\begin{align}\label{eq:gdecayAS}
g_{ij}&=(g_{m})_{ij}+{O}_{2}(r^{-\tau}),\\\label{eq:KdecayAS}
K_{ij}&=\phantom{(g_{m})_{ij}+}\;{O}_{1}(r^{-(\tau+1)})
\end{align}
with $g_{m}$ as in Section~\ref{sec:noRT} for some $\tau>1$ (automatically implying \eqref{eq:Hamint}, \eqref{eq:momint} and indeed the weak Regge--Teitelboim conditions). This condition is known as being \emph{asymptotic to (spatial) Schwarzschild} or \emph{asymptotically Schwarzschildean}; it is more typically written in so-called \emph{isotropic coordinates} but this is not central here. Of course, in line with our results presented in Section~\ref{sec:noRT}, not all asymptotically Euclidean initial data sets are asymptotically Schwarzschildean. 

In view of the question whether asymptotic flatness can be characterized geometrically to be addressed in Section~\ref{sec:geometric}, one can ask whether/how asymptotic spherical symmetry in the above sense can be characterized geometrically. An important step towards this goal has been achieved in this project by Avalos~\cite{avalos2024sobolevregularitycompactified3manifolds}:

Loosely put, for a time-symmetric asymptotically Euclidean initial data set $(M,g,K=0)$ (with $\tau<1$), he identifies the following sufficient conditions for being asymptotically Schwarzschildean: 
\begin{itemize}
\item the scalar curvature $\scal_{g}$ satisfies $r^{\eps}\scal_{g}\in L^{p}(\R^{3}\setminus\overline{B})$ for some $\eps>0$ and $p>3$, and
\item the Cotton tensor $C_{g}$ of $g$ satisfies $r^{-\sigma}(C_{g})_{ijk}\in L^{q_{\sigma}}(\R^{3}\setminus\overline{B})$ for some $-3<\sigma<-1$ and $q_{\sigma}\definedas\frac{3}{6+\sigma}$
\end{itemize}
(with respect to the given asymptotic coordinates)\footnote{It might be interesting to notice that there is a parallel concerning regularity. In \cite{AvalosCogoAbrego}, it has been recognized how a combined control on the regularity of $\scal_g$ and $C_g$ controls the optimal regularity element of the conformal class $[g]$ for $g\in W^{2,q}(M)$, $q>n$, and $M$ closed. Thus, the regularity-decay controls of $(\scal_g,C_g)$ in \cite{avalos2024sobolevregularitycompactified3manifolds,AvalosCogoAbrego}  can be presented as a conformal parallel of the fact that regularity of $\mathrm{Ric}_g$ controls local regularity due to \cite{DeTurckKazdan} while it controls asymptotic decay on asymptotically Euclidean manifolds by \cite{Bartnik86}.}. Furthermore, Avalos~\cite{avalos2024sobolevregularitycompactified3manifolds} asserts that these conditions are also sufficient for $\vec{Z}_{\text{B\'ORT}}$ to converge. The proof relies on a very intricate geometric construction of coordinates relying on suitably regular conformal compactifications as well as on the use of Green's functions.

As a spin-off of this project, Avalos and Cederbaum~\cite{AvalosCederbaum} are studying how ``sharp'' the above sufficient conditions are for being asymptotically Schwarzschildean. In particular, they find examples of asymptotically Euclidean initial data of order $\tau=1$ which satisfy the scalar curvature assumption but barely fail the Cotton tensor assumption; these examples do not possess any coordinates in which they are asymptotically Schwarzschildean.


\section{Center of Mass and Geometric Foliations}
\label{sec:foliation}
Taking a more geometric approach to defining the center of mass, Huisken and Yau \cite{Huisken-Yau-1996} proposed to use \emph{asymptotic constant mean curvature (CMC-) foliations} to take the role of the geometric center of mass
of an asymptotically flat Riemannian manifold or initial data set. They showed that in every $3$-dimensional Riemannian manifold $(M,g)$ which is asymptotic to a Riemannian 
Schwarzschild manifold of positive mass (a condition much stronger than \eqref{eq:gdecay}, namely requesting \eqref{eq:gdecayAS} with $\tau=2$, and with four derivatives instead of just two), there exists a smooth foliation
$\{\Sigma_s^{\text{CMC}}\mid \sigma \in(\sigma_0,\infty)\}$ of the asymptotic end
by stable surfaces of constant mean curvature $H(\Sigma_\sigma^{\text{CMC}}) = \frac{2}{\sigma}$. Furthermore, they proved conditional uniqueness for the leaves of the foliation.

Subsequent works by Ye \cite{Ye}, Huang \cite{Huang-2010}, Eichmair
and Metzger \cite{Eichmair-Metzger-2013}, and Nerz
\cite{Nerz-CMC-Foliation-2015} etc. were able to
considerably relax the decay assumptions for the initial data under
which these CMC-foliations exist and also give various uniqueness
results. To be precise, the most general existence result in dimension
$n=3$ is by Nerz \cite[Theorem 5.1 and Theorem
5.3]{Nerz-CMC-Foliation-2015}. It holds under the
asymptotic assumptions~\eqref{eq:gdecay}, \eqref{eq:Hamint} for $\tau>\tfrac{1}{2}$ if the ADM-energy~\eqref{def:E} of $(M,g)$ is non-vanishing.

The previous results on the existence of CMC-surfaces are almost
exclusively in $3$ dimensions, except \cite{Eichmair-Metzger-2013}
in the case of asymptotically Schwarzschildean manifolds. For general
dimension $n\geq3$ and general asymptotics
$(\Phi_*g)_{ij} = \delta_{ij} + O_3(r^{-\tau})$ for
$\tau\in(\frac{n-2}2,n-2]$, Eichmair and K\"orber 
\cite{eichmair2024foliationsasymptoticallyflat3manifolds} give a
conceptually new proof that makes use of a Lyapunov--Schmidt
reduction.

Tenan and Sinestrari~\cite{TenanSinestrari} picked up the original strategy of proof by Huisken and Yau~\cite{Huisken-Yau-1996} and reproved Nerz's and K\"orber--Eichmair's ($n=3$) results via volume-preserving inverse mean curvature flow, carefully adjusting all estimates to the asymptotic flatness condition $\tau>\frac{1}{2}$ (under a technical assumption related to the weak Regge--Teitelboim conditions).

The crucial feature of all these existence proofs is the fact that all spheres have constant mean
curvature in Euclidean space. Hence uniqueness only holds up to Euclidean motions. In an
asymptotically flat manifold, this causes the linearization of the
mean curvature operator to have an approximate kernel, which
complicates the use of the inverse function theorem. If the ADM-energy~\eqref{def:E}
does not vanish, a careful analysis shows that this kernel in
fact disappears, provided the decay conditions are strong
enough. This is roughly the base of Nerz's proof.

In contrast, Eichmair and K\"orber first generate a family of surfaces
$\Sigma_r(a)$ which arise from spheres $\mathbb{S}_r(a)$ of radius $r$ centered
at some point $a$ by solving for a graphical perturbation of $\mathbb{S}_r(a)$
such that the $\Sigma_r(a)$ are CMC up to the approximate kernel. Then
they carefully analyze the area function $G\colon a\mapsto |\Sigma_r(a)|$ and
show that a critical point $a_r$ exists and that $\Sigma_r(a_r)$ has
constant mean curvature.

Considering uniqueness of the foliation, Qing and Tian \cite{QT} showed
that in an asymptotically Schwarzschildean manifold with positive mass,
stable surfaces of constant mean curvature enclosing a large enough
compact set are part of the foliation $\{\Sigma_\sigma^{\text{CMC}}\}$. This was
generalized by Ma \cite{Ma} to metrics of the form
$(\Phi_*g)_{ij} = \delta_{ij} + O_4(r^{-1})$ and positive energy.

For some time, only conditional uniqueness results were available when
relaxing the decay rates or the number of derivatives required in
these uniqueness theorems. To state these, for a surface
$\Sigma\subset M\setminus C$ denote by $\rho(\sigma)\definedas  \min \{\Phi(p) \mid p\in \Sigma\}$ the distance of
$\Sigma$ to the origin with respect to the diffeomorphism $\Phi$. The uniqueness theorem by Nerz \cite[Theorem
5.3]{Nerz-CMC-Foliation-2015} only requires~\eqref{eq:gdecay} with
non-vanishing energy.
 It states that for every constant $c$ there is a
constant $\sigma_0$ such that if $\Sigma$ is a stable surface with
constant mean curvature $\frac{2}{\sigma}$ for $\sigma>\sigma_0$, if
$\Sigma$ encloses $C$, and if
\begin{equation}
  \label{eq:balancing}
  \frac{2}{H} <c\rho(\Sigma)
\end{equation}
then $\Sigma$ is part of the foliation from the existence theorem,
that is, $\Sigma=\Sigma_\sigma^{\text{CMC}}$.

If one assumes suitably strong Regge--Teitelboim conditions for the metric (and its
derivatives of order up to five), then condition~\eqref{eq:balancing}
can be relaxed to
\begin{equation}
  \label{eq:balancing2}
  \frac{2}{H} < \rho(\Sigma)^s
  \quad\text{for}\quad
  s\in\left(1,\frac{4+2\tau}{5-\tau}\right).
\end{equation}
This is the uniqueness theorem by Huang \cite{Huang-2010}, and
here $\sigma_0$ depends on $s$. 

Eichmair and K\"orber removed the conditions of
type~\eqref{eq:balancing} and ~\eqref{eq:balancing2} and gave a
complete answer to the uniqueness question in form of an unconditional
uniqueness theorem
\cite[Theorem 12]{eichmair2024foliationsasymptoticallyflat3manifolds}
in dimension $n=3$ assuming non-negative scalar curvature. Their
Theorem 14 also improves the conditional uniqueness theorems in the
case when the scalar curvature is allowed to change sign to the range
$s\in\left(1,1+\frac{3(2\tau-1)}{4(1-\tau)}\right)$ without assuming
Regge--Teitelboim conditions.

\subsection{A coordinate center for the CMC-foliation}
To attach a center of mass coordinate vector $\vec{Z}_{\text{CMC}}$ to the asymptotic CMC-foliation, Huisken and Yau
\cite{Huisken-Yau-1996} computed the barycenters $z(\sigma)$ of the leaves
$\{\Sigma_\sigma\}$ in asymptotic coordinates $x$, that is,
\begin{equation*}
  \vec{z}(\sigma) \definedas  \frac1{|\Sigma_\sigma^{\text{CMC}}|_{\delta}}\int_{\Sigma_\sigma^{\text{CMC}}} x \,dA_\delta
\end{equation*}
and then defined the limit $\vec{Z}_\text{CMC} \definedas  \lim_{\sigma\to\infty}
\vec{z}(\sigma)$ to be the \emph{CMC-center of mass} --- provided it exists. 

Assuming non-vanishing energy as well as the strong Regge--Teitelboim
conditions on the metric, its derivatives, and the curvature,
Huang \cite{Huang-2010} was able to show that the limit $\vec{Z}_\text{CMC}$
exists and is equal to $\vec{Z}_{\text{B\'ORT}}$. Recall that the latter
converges due to the strong Regge--Teitelboim
conditions. Nerz~\cite{Nerz-CMC-Foliation-2015} refined this by
assuming only the weak Regge--Teitelboim conditions which do not
guarantee the convergence of $\vec{Z}_{\text{B\'ORT}}$. Instead he
showed that $\vec{Z}_\text{CMC}$ converges if an only if $\vec{Z}_{\text{B\'ORT}}$ converges and then they coincide. The higher dimensional case is again treated in Eichmair and K\"orber \cite[Theorem
10]{eichmair2024foliationsasymptoticallyflat3manifolds}.

As noted by Cederbaum and Nerz \cite{Cederbaum-Nerz-2015}, the Regge--Teitelboim conditions are crucial for this convergence to hold and even in
physically relevant data, such as (graphical) initial data sets in the
Schwarzschild spacetime --- see Section~\ref{sec:noRT} ---, the barycenters of the CMC-leaves
$\{\Sigma_\sigma^{\text{CMC}}\}$ can oscillate and fail to converge. This is not a
genuine phenomenon in relativity though, as such phenomena can be
observed already in Newtonian gravity, as described in
\cite[Section 3]{Cederbaum-Nerz-2015} and
\cite{Cederbaum2023}.

\begin{remark}[Asymptotic equivariance under Euclidean motions]
Note that the CMC-center of mass $\vec{Z}_\text{CMC}$ naturally has the desired covariance
properties under asymptotic rotations and translations of the
coordinate chart, so that in view of the above established equality,
the B\'ORT-center of mass $\vec{Z}_\text{B\'ORT}$ also transforms accordingly -- a fact already established by Chru\'sciel~\cite{Chrusciel} under the assumption of strong Regge--Teitelboim conditions.
\end{remark}

\begin{remark}[Asymptotic Poincar\'e equivariance]\label{rem:boostequivariance}
In a similar spirit, Chru\'sciel~\cite{Chrusciel} also establishes that this center of mass transforms as desired under asymptotic Poincar\'e transformations once one assumes strong Regge--Teitelboim conditions. Again, this results is established with respect to one fixed asymptotic coordinate system, a technical point that might actually matter, see forthcoming work by Senthil Velu~\cite{Sara} and by Cederbaum~\cite{Csuper} as well as Section~\ref{sec:ang}.
\end{remark}

\subsection{Including the spacetime perspective: the STCMC-foliation}\label{sec:STCMC}
From a physics perspective, it might appear strange that the second fundamental form $K$ of an initial data set $(M,g,K)$ does not appear in either of the above definitions of center of mass, see Szabados~\cite{Szabados}. Intuitively speaking and very much simplified, this is because $K$ can be regarded as capturing the infinitesimal  velocity of the system as measured by the asymptotic ``observers'' -- which should matter for locating the center of mass in view of special relativity. Indeed, this manifests itself when studying asymptotic (coordinate) boosts, see also Section~\ref{sec:ang}.

To address this issue, Cederbaum and Sakovich~\cite{Cederbaum-Sakovich:2018a} introduced a new foliation which
requires the leaves $\SigST$ to have constant \emph{spacetime mean
curvature (STCMC)}
\begin{equation}
  \label{eq:STCMC}
  \mathcal{H}(\Sigma) \definedas  \sqrt{H^2 - (\tr_\Sigma K)^2 } = \tfrac2\sigma. 
\end{equation}
Here, $\tr_\Sigma K$ is the $2$-dimensional trace of $K$ restricted to
the tangent space of~$\Sigma$.  The spacetime mean curvature is a spacetime
quantity in the sense that $\mathcal{H}(\Sigma)$ is the Lorentzian length of
the mean curvature vector of the embedding of $\Sigma$ into the
spacetime constructed from $(M,g,K)$ --- provided the mean curvature vector is spacelike --- and thus does not depend on the
choice of spatial slice $M$ passing through $\Sigma$. 

Existence and uniqueness of the STCMC-foliation for asymptotically Euclidean initial data sets for any $\tau>\frac{1}{2}$ with $E\neq0$ are asserted in \cite{Cederbaum-Sakovich:2018a} via a method of continuity argument based\footnote{We would like to remark that \cite{Nerzevo,Nerz-CMC-Foliation-2015,NerzCE} contain inaccuracies and small gaps which were addressed and fixed in \cite{Cederbaum-Sakovich:2018a}, see also forthcoming work by Olivia Vi\v{c}\'{a}nek Mart\'{i}nez~\cite{VDiss}.} on the result of Nerz~\cite{Nerz-CMC-Foliation-2015} and on the construction of asymptotic constant expansion foliations by Metzger~\cite{Metzger-2007} later refined by Nerz~\cite{NerzCE}. The uniqueness assertion \cite[Theorem 4]{Cederbaum-Sakovich:2018a} is a natural extension of the one by Nerz, see the beginning of Section~\ref{sec:foliation}. Another proof of this result was given by Tenan~\cite{Tenan}, generalizing \cite{TenanSinestrari} to the initial data context by studying volume-preserving spacetime mean curvature flow.

Once the STCMC-foliation $\{\SigST\}$ is constructed, one can
again study the barycenters $\vec{z}(\SigST)$ in an asymptotic
coordinate chart $x$ and study their limit 
\begin{equation*}
  \vec{Z}_\text{STCMC} \definedas  \lim_{\sigma\to\infty} \vec{z}(\SigST)
\end{equation*}
called the \emph{STCMC-center of mass} provided it converges. As it turns out \cite[Theorem
5]{Cederbaum-Sakovich:2018a}, given the decay
conditions~\eqref{eq:gdecay}--\eqref{eq:momint} and assuming in addition the strong
Regge--Teitelboim condition on $g$, combined with $K=O(r^{-2})$ as $r\to\infty$, $\vec{Z}_\text{STCMC}$ converges if and only if the
following surface integral involving the \emph{conjugate momentum tensor} 
\begin{equation}\label{def:pi}
\pi \definedas  (\tr K) g - K
\end{equation}
and the ADM-energy $E$ converges:
\begin{equation}
  \label{eq:correction}
  Z_0^i \definedas  \lim_{r\to\infty} Z_0^i(\mathbb{S}^{2}_r)
  \qtext{with}
  Z_0^i(\Sigma) \definedas   \frac{1}{32\pi E} \int_{\Sigma}
  \tfrac{x^i}{r}\, [\pi \left(\tfrac{x}{r},
    \tfrac{x}{r}\right)]^{2}
  \,d A_\delta.
\end{equation}
The vector $\vec{Z}_0$ then is a correction to the B\'ORT/CMC-center of mass in
the sense that 
\begin{equation}
  \label{eq:stcmc}
  \vec{Z}_\text{STCMC} = \vec{Z}_\text{B\'ORT} + \vec{Z}_0=\vec{Z}_{\text{CMC}}+\vec{Z}_{0}.
\end{equation}
We would like to point out that the correction term $\vec{Z}_{0}$ vanishes entirely under the full strong Regge--Teitelboim conditions as can be seen by parity arguments.

\begin{remark}[Time-evolution]\label{rem:evo}
The STCMC-center of mass evolves as the position of a point particle in general relativity, that is, it satisfies~\cite[Theorem
6]{Cederbaum-Sakovich:2018a} 
\begin{equation*}
  \frac{d}{dt} \vec{Z}_\text{STCMC} = \frac {\vec{P}}{E}.
\end{equation*}
This statement also holds asymptotically in a leaf-wise sense without assuming Regge--Teitelboim conditions and thus does not require $\vec{Z}_{\text{STCMC}}$ to converge. Interestingly, the same evolution law also holds for the CMC-foliation as was established by Nerz~\cite{Nerzevo}; this is consistent as one computes that indeed $\frac{d}{dt}\vec{Z}_{0}=0$. See also forthcoming work by Senthil Velu~\cite{Saraevo} who removes the technical assumption $K=O(r^{-2})$ as $r\to\infty$ in \cite{Nerzevo} and \cite[Theorem
6]{Cederbaum-Sakovich:2018a}.
\end{remark}

\begin{remark}[Asymptotic boost equivariance]
As the STCMC- and B\'ORT/CMC-centers of mass coincide under strong Regge--Teitelboim conditions, Remark~\ref{rem:boostequivariance} on page~\pageref{rem:boostequivariance} also applies to the STCMC-center of mass. We have indications that the STCMC-center of mass transformation behavior under asymptotic coordinate boosts fits better with the expectations from special relativity than that of the B\'ORT/CMC-center of mass, see forthcoming work by Senthil Velu~\cite{Sara} and Section~\ref{sec:ang}.
\end{remark}

A first indication of the conjecture that the STCMC-foliation captures the center of mass better than the CMC-foliation comes from the explicit example by Cederbaum and Nerz~\cite{Cederbaum-Nerz-2015}, as was analyzed by Cederbaum and Sakovich~\cite[Section 9]{Cederbaum-Sakovich:2018a}. Recalling the graphical setup in the Schwarzschild spacetime from Section~\ref{sec:noRT}, this corresponds to choosing the graph function 
\begin{equation*}
  f \colon  M_m \to \R \colon x \mapsto \sin(\ln r) + \tfrac{\< x, a \>}{r}
\end{equation*}
for some fixed $0\neq a\in\R^{3}$ (corresponding to $\beta=0$, $\gamma=0$ in Section~\ref{sec:noRT}). The graphical initial data set $(M_{f},g_{f},K_{f})$ is asymptotically Schwarzschildean of order $\tau=2$ and hence satisfies the assumptions of  \cite[Theorem
5]{Cederbaum-Sakovich:2018a}. In particular, $\vec{Z}_{\text{B\'ORT}}=\vec{Z}_{\text{CMC}}$ holds. An explicit computation based on \eqref{eq:BORTcenter} shows that
\begin{equation*}
  \vec{Z}_\text{B\'ORT}(\mathbb{S}^{2}_r) = \tfrac{\cos(\ln r)}{3m}\,a + O(r^{-1})
\end{equation*}
as $r\to\infty$. In particular, this does not converge as $r\to\infty$ and hence $\vec{Z}_{\text{CMC}}$ does not converge either. In contrast, one finds 
\begin{equation*}
  \vec{Z}_\text{B\'ORT}(\mathbb{S}^{2}_r) + \vec{Z}_0(\mathbb{S}^{2}_r) = O(r^{-1})  
\end{equation*}
as $r\to\infty$ so that $\vec{Z}_{\text{STCMC}}$ indeed \emph{does} converge to the coordinate origin as should be expected from the spherical symmetry of the Schwarzschild spacetime.

\subsection{Related results on prescribed mean curvature foliations}
CMC-, STCMC-, and related foliations have been studied in other contexts as well. While it would lead too far to highlight all those contexts, we would like to mention those which are most intimately related with this project. 

Staying in the context of initial data sets $(M,g,K)$ and in particular Riemannian manifolds $(M,g)$, an analytically different regime where it is interesting to study
foliations of prescribed mean curvature is the local case. Here, we ask
the question whether a foliation by e.g. constant mean curvature
spheres exists near a point $p\in M$, where $(M,g)$ is an arbitrary
Riemannian manifold. The first result in this direction is the
construction by Ye~\cite{Ye1991} who showed that if $p$ is a
non-degenerate critical point of the scalar curvature $\scal_g$ of $g$,
then there exists a neighborhood $U$ of $p$ such that
$U\setminus\{p\}$ is foliated by CMC-surfaces. The condition that $p$ be a critical point is actually a
necessary condition for such a foliation to exist. 

Ye's method uses a kind of Lyapunov--Schmidt reduction. Constructing foliations with different prescribed mean curvature leaves on the basis of similar reductions yield other necessary conditions which have to be satisfied at the
concentration point $p$ of the leaves.

As part of this project, Metzger and Pe\~nuela studied the local
STCMC-foliation, where leaves satisfy the STCMC-equation
\eqref{eq:STCMC} as well as the local foliation by surfaces of constant expansion
(CE) in \cite{MetzgerPenuela23}, where the leaves satisfy the
equation
\begin{equation}
  \label{eq:CEfoliation}
  H \pm \tr_\Sigma K = \tfrac{2}{\sigma}
\end{equation}
where one is interested in the case of small
$\sigma\in(0,\sigma_0)$. Note that the choice of sign in
\eqref{eq:CEfoliation} gives rise to two different foliations of
surfaces with constant expansion. 

The existence of the local STCMC-foliation is governed by the
\emph{local STCMC $1$-form} \cite[Definition 1.4]{MetzgerPenuela23}
which is defined (in dimension $n=3$) by
\begin{equation*}
    A_\text{ST}  =\tfrac{3}{2} \nabla \scal_g + \tfrac{1}{8}\Big[
    \nabla |K|^2  + \tfrac{41}{2} \nabla (\tr K)^2
    +4 \diver ( K^2 )  - 14 \diver \left(\tr K\cdot K\right)\Big].
\end{equation*}
It replaces the gradient of the scalar curvature in the original
result by Ye.

For a full statement of the existence result, we refer the reader to \cite[Section
1.1]{MetzgerPenuela23}, but roughly speaking \cite[Theorem
3.2]{MetzgerPenuela23} and \cite[Theorem 5.4]{MetzgerPenuela23} give
that a unique foliation by STCMC-surfaces exists near $p\in M$ provided that
\begin{enum}
\item $A_\text{ST}(p)=0$,
\item $\nabla A_\text{ST}(p)$ is invertible, and
\item $|K|$ and $|\nabla K|$ are sufficiently small near $p$.
\end{enum}
It remains open to give a  geometric and/or
physical interpretation of $A_\text{ST}$.

Analytically the case of constant expansion foliations is of a
different nature, since the term $\tr_\Sigma K$ in
\eqref{eq:CEfoliation} enters at a different order. It thus
turns out that only $K$ determines the corresponding one-form which
determines the existence of a CE-foliation near $p\in M$
(cf. \cite[Definition 1.5]{MetzgerPenuela23}):
\begin{equation*}
  A_\text{CE}  =\tfrac{5}{6} \nabla \tr K -2 \diver K. 
\end{equation*}

Again, the full statement for the existence and uniqueness of a
CE-foliation near $p$ which is proved in \cite[Theorem
4.3]{MetzgerPenuela23} and \cite[Theorem 5.5]{MetzgerPenuela23} is too
technical, but a unique foliation by CE-surfaces exists near $p\in M$
provided that
\begin{enum}
\item $A_\text{CE}(p)=0$, $K(p) =0$,
\item $\nabla A_\text{CE}(p)$ is invertible,
\item $|\nabla R_g|$, $|k|$ and $|\nabla k|$ are sufficiently small near $p$.
\end{enum}

Moving away from initial data sets, Kr\"oncke and Wolff~\cite{kroenckewolff} construct asymptotic foliations by STCMC-surfaces  in an asymptotically Schwarzschildean lightcone. Motivated by the fundamental work of Huisken and Yau \cite{Huisken-Yau-1996}, they employ a modified mean curvature flow along the null hypersurface. While their result relies on strong asymptotic assumptions as in the initial work of Huisken and Yau, they develop all the necessary tools and crucial geometric insights to extend this result to more general asymptotics in the null setting. As in the case of initial data sets, they propose a geometric notion of center of mass motivated by special relativity. Their results so far suggest that this geometric definition might yield an approach to address the super-translation ambiguity of the classical notion of center of mass at null infinity, and they are currently investigating the properties of their center of mass definition and the uniqueness of their foliation in an ongoing work. Extending these considerations of center of mass to the null setting is potentially of high physical relevance, as any information such as radiation reaching us from a far away celestial object travels in spacetime along such a null hypersurface.


\section{Constructing geometric coordinates from geometric foliations}
\label{sec:geometric}

The reciprocal question to the existence results discussed in
Section~\ref{sec:foliation}, namely to construct asymptotically Euclidean coordinates from an existing geometric foliation was first studied by Nerz in \cite{Nerz-AF-2015}. He argues that if a Riemannian $3$-manifold admits a CMC-foliation which satisfies certain geometric
conditions then it is asymptotically flat. This gives a geometric
characterization of asymptotic flatness without reference to an
asymptotic coordinate chart similar to Bando, Kasue, and Nakajima
\cite{BKN1989}, who relied on pointwise curvature decay and volume
growth assumptions.

More precisely, Nerz explains why a $3$-dimensional Riemannian manifold is
asymptotically Euclidean if it possesses a weak CMC-cover satisfying a set
of geometric and analytic conditions. These include local uniqueness,
controlled instability, and curvature decay estimates. His proof relies on the explicit construction of asymptotic coordinates\footnote{We would like to remark that \cite{Nerz-AF-2015} contains some unfortunate typos and oversights, see forthcoming work by Piubello and Vi\v{c}\'{a}nek Mart\'{i}nez~\cite{PVM} and by Vi\v{c}\'{a}nek Mart\'{i}nez~\cite{VDiss}.}.

One part of this project, pursued by Piubello and
Vi\v{c}\'anek Mart\'{i}nez~\cite{PVM}, generalizes these ideas to the
more general setting of an initial data set $(M, g, K)$, and
to foliations by STCMC-surfaces, as discussed in Section~\ref{sec:STCMC}. They consider a family of $2$-dimensional STCMC-spheres that weakly foliate $M$. In view of the discussion in Section \ref{sec:STCMC}, one can expect that the coordinates constructed from these surfaces are more suitable for analyzing quantities related to the full spacetime and their time evolution and transformation properties under asymptotic boosts.

The strategy to construct asymptotically Euclidean coordinates is to
adapt the construction suggested by Nerz to this new setting. That is,
consider a family $\{\Sigma_\sigma^{\text{STCMC}}\}$ of STCMC-surfaces labeled by
their \emph{spacetime mean curvature radius $\sigma$}, so that each leaf has constant spacetime
mean curvature  $\mathcal{H}(\Sigma_{\sigma}^{\text{STCMC}}) =\frac{2}{\sigma}$. The construction relies
on several assumptions about this family of surfaces, mostly geometric
bounds:
\begin{enum}
\item a uniform bound on the Hawking mass of the leaves $\Sigma^{\text{STCMC}}_\sigma$,
\item control on the decay of the Ricci curvature, the second
  fundamental form $K$, and the energy and momentum densities, 
\item local uniqueness of the foliation,
\item the leaves are pairwise  disjoint, and their union covers the
  complement of a compact set of $M$,
\item each of the $\Sigma_\sigma^{\text{STCMC}}$ is such that smallest eigenvalue of
  the stability operator is bounded from below in a controlled way.
\end{enum}
The last assumption is weaker than requiring stability and allows to
study settings with negative mass. 

The proof begins by establishing a priori estimates on each leaf
$\Sigma_\sigma^{\text{STCMC}}$. The main step is to show that the trace-free part of
the second fundamental form is small in $L^4$, which implies that the
surfaces are asymptotically almost umbilic. Then an argument \`{a} la
De Lellis--M\"uller~\cite{DeLellisMueller}, used by Nerz in the Riemannian setting, applies and
gives a conformal parametrization of each leaf with a 
close to constant conformal factor. As a result, the geometry of the STCMC-foliation can be
estimated to be close to that of centered round spheres.

Furthermore, these estimates imply that the STCMC-stability operator has trivial kernel
on each leaf and thus, by the Fredholm
alternative, its adjoint also has trivial kernel and the operator is
invertible. Applying the inverse function theorem to the spacetime mean curvature
map gives a diffeomorphism
\begin{equation*}
  \Psi\colon (\sigma_0, \infty) \times \mathbb{S}^2 \to M, \quad \Psi(\sigma, \mathbb{S}^{2}) = \Sigma_\sigma^{\text{STCMC}}.
\end{equation*}
A delicate analysis performed in~\cite{PVM} leads to asymptotic control on the lapse function. 
The geometry of the STCMC-foliation is now known to be regular enough to compare the
Laplacian on each leaf $ \Sigma_\sigma^{\text{STCMC}}$ to that of the standard round
sphere. In particular, on each leaf, this Laplace operator has three distinct
eigenfunctions $ f_1^\sigma $, $ f_2^\sigma $, and $ f_3^\sigma $,
corresponding to the linear eigenfunctions on $ \mathbb{S}^2 $. These
are carefully extended smoothly along the foliation using the map
$\Psi$ via a suitable minimal rotation argument.

The asymptotic coordinates can then be defined in terms of these
eigenfunctions, appropriately rescaled and augmented by a correction
term derived from the lapse function of the STCMC-foliation.  A careful
analysis shows that these coordinates are well-defined and that the
components of the metric expressed in this system have the expected
asymptotically Euclidean decay \eqref{eq:gdecay}, both along the leaves and in the
radial direction. This completes the construction.

Combining this construction with the existence and uniqueness result
for STCMC-foliations~\cite{Cederbaum-Sakovich:2018a,Tenan}, this gives a fully geometric
characterization of asymptotically Euclidean initial data sets, both
in the pointwise and the weak sense.

A natural application of this coordinate construction is the study of the center of mass in the constructed coordinates. In particular, it is very interesting to study the convergence of the surface
integral in the B\'ORT-center of mass definition~\eqref{eq:BORTcenter} and of its correction using the
term~\eqref{eq:stcmc}. This might yield a resolution of \cite[Conjecture 1]{Cederbaum-Sakovich:2018a} roughly conjecturing that there should be a geometric condition on asymptotic coordinates $x$ ensuring that $\vec{Z}_{\text{STCMC}}$ will converge for asymptotic coordinates $x$ if $x\rho\in L^{1}(\R^{3}\setminus\overline{B})$. While Avalos' results~\cite{avalos2024sobolevregularitycompactified3manifolds} described in Section~\ref{sec:AS} provide some indication in favor of this conjecture, the coordinates constructed by Piubello and Vi\v{c}\'anek Mart\'{i}nez~\cite{PVM} may be a promising candidate for being the conjectured coordinates $x$ in full generality.


\section{Boosts and angular momentum: problems and first ideas}\label{sec:ang}
In Section~\ref{sec:CoM}, we have introduced the B\'ORT-expression that arises as the definition of center of mass by the Hamiltonian approach taken by Regge and Teitelboim~\cite{Regge--Teitelboim-1974} and Beig and \'O Murchadha~\cite{Beig-OMurchadha-1987} and coinciding with the definition given by Michel~\cite{Michel}. Regge and Teitelboim~\cite{Regge--Teitelboim-1974} also derive a (Hamiltonian) definition of angular momentum $\vec{J}_{\text{RT}}$ for asymptotically Euclidean initial data $(M,g,K)$ with respect to asymptotic coordinates $x$. The components of $\vec{J}_{\text{RT}}$ are given by
\begin{equation}\label{eq:JRT}
J^k \definedas  \frac{1}{8\pi}\lim_{r \rightarrow \infty}\int_{\mathbb{S}^2_r}\pi(Y^{(k)},\tfrac{x}{r}) \; dA_{\delta},
\end{equation}
$k=1,2,3$, where $Y^{(3)}=x^{1}\partial_{x^{2}}-x^{2}\partial_{x^{1}}$ with $Y^{(1)}$, $Y^{(2)}$ defined accordingly via cyclic permutations and where $\pi$ is the conjugate momentum tensor defined in~\eqref{def:pi}. Just as for $\vec{Z}_{\text{B\'ORT}}$, $\vec{J}_{\text{RT}}$ does not converge in general but is well-known to converge under strong Regge--Teitelboim conditions. Similarly, Michel~\cite{Michel} also provides a definition of angular momentum $\vec{J}_{\text{M}}$ with components given by
\begin{equation}\label{eq:JM}
J^k \definedas  \frac{1}{8\pi}\lim_{r \rightarrow \infty}\int_{\mathbb{S}^2_r}\overline{\pi}(Y^{(k)},\tfrac{x}{r}) \; dA_{\delta},
\end{equation}
where $\overline{\pi}_{ij}\definedas K_{ij}-(\tr_{\delta}K)\delta_{ij}$ \emph{depends} on the choice of asymptotic coordinates~$x$ in contrast to the definition of $\pi$ in \eqref{def:pi}. Upon first inspection and according to a comment in~\cite{Michel}, it might appear that $\vec{J}_{\text{RT}}$ and $\vec{J}_{\text{M}}$ coincide (whenever they converge) but this is not the case as is observed by Cederbaum~\cite{Csuper}, where explicit examples of initial data sets are given (as graphical slices in the Schwarzschild spacetime) for which both expressions converge but not to the same limit. Some of these examples were first studied from a different perspective by Chen, Huang, Wang, and Yau~\cite{Chen2016}. However, both expressions \emph{do} give the same result when strong Regge--Teitelboim conditions are assumed.

We are now in a position to make our previous remarks on asymptotic Poincar\'e equivariance a little more precise. In fact, just as one combines energy and linear momentum into an energy-momentum $4$-vector which is then a spacetime invariant, one combines angular momentum and the \emph{center of mass charge} $t\vec{P}-E\vec{Z}_{\text{B\'ORT}}$ into an antisymmetric \emph{angular momentum $4$-tensor} of type $(0,2)$. This spacetime angular momentum tensor is then expected to be spacetime equivariant in the sense that it is invariant under the asymptotic Lorentz group and transforms equivariantly under asymptotic translations. This has been established in a rigorous sense by Chru\'sciel~\cite{Chrusciel}, with respect to a fixed asymptotic chart and assuming strong Regge--Teitelboim conditions. Chru\'sciel's results hence apply to all the notions of center of mass discussed above (entering into the center of mass charge). A central ingredient in his proof is the representation of the corresponding charges via
super-potentials and a Stokes' theorem argument over a timelike cylindrical region. 

Leaving the realm of the strong Regge--Teitelboim conditions and inspecting the definitions of RT-/M-angular momentum and B\'ORT-center of mass more closely, Cederbaum~\cite{Csuper} asserted that there is a subtle fundamental problem underlying these definitions. This problem is due to the fact that the coordinate dependent ``asymptotic Killing vector fields'' inserted into the Hamiltonian in the Hamiltonian approach do not generically asymptotically satisfy the Killing equations to a satisfactory order in the absence of strong Regge--Teitelboim conditions; this is related to the super-translations mentioned in Remark~\ref{rem:supertrans} on page~\pageref{rem:supertrans} (or rather a spacetime version thereof), see \cite{Csuper} for more information. Moreover, the same problem occurs in Michel's approach where he uses asymptotic Killing initial data to be understood as lapse and shift of the same ``asymptotic Killing vector fields''. This leads us to suspect that the definitions of angular momentum and center of mass might need to be adjusted in the absence of strong Regge--Teitelboim conditions. 

It is not at all clear whether the geometric approach via the STCMC-foliation described in Section~\ref{sec:STCMC} remedies this problem for the definition of center of mass; this is currently researched by the first author as a continuation of this project. It is even less clear how to remedy the problem by geometric means for the definition of angular momentum. First ideas in the direction of the second question are being pursued by Vi\v{c}\'anek Mart\'{i}nez as a continuation of this project.

Coming back to asymptotic Poincar\'e equivariance, recall that we have pointed out that Chru\'sciel's result are obtained in a fixed asymptotic coordinate system. Specifically, this means that asymptotic Poincar\'e transformations are interpreted as Poincar\'e transformations performed on the asymptotic coordinates (not allowing for any lower order terms). In the example discussed in \eqref{eq:transfo}, this would mean that one only allows lower order terms that are exact translations, see Remark~\ref{rem:supertrans}~on page~\pageref{rem:supertrans}. 

Investigating the existing definitions of center of mass further, it is thus very relevant to understand how they transform under more general asymptotic coordinate transformations which are only asymptotic to boosts, in line with the evolution results described in Remark~\ref{rem:evo} on page~\pageref{rem:evo}, and in the absence of Regge--Teitelboim conditions. This is pursued in forthcoming work by Senthil Velu~\cite{Sara} as part of this project and asserts the following results regarding the transformation behavior of the B\'ORT-/CMC- and STCMC-center of mass under infinitesimal boosts: 

One of her results relies on an observation by Chru\'sciel~\cite{Chruscielangular} stating that $\vec{J}_{\text{RT}}$ converges for asymptotically Euclidean initial data sets $(M,g,K)$ for any asymptotic coordinates $x$ in which 
\begin{align}
g_{ij}&=\delta_{ij}+{O}_{2}(r^{-p}),\\
K_{ij}&=\phantom{\delta_{ij}+}\;{O}_{1}(r^{-q})
\end{align}
hold for some $p+q>3$ for which $p>\frac{1}{2}$, $q>\frac{3}{2}$. Then, for any fixed unit vector $n\in(\R^{3},\delta)$ --- to be thought of as specifying a boost direction ---, the STCMC-center of mass $\vec{Z}_{\text{STCMC}}$ changes upon infinitesimal boosts by the transformation law
\begin{equation}
\frac{d}{d\omega} Z_{\text{STCMC}}^{\,i}
= \frac{1}{E}\,\epsilon_{ijk} \: n^j J_{RT}^k,
\end{equation}
where $\epsilon_{ijk}$ denotes the totally antisymmetric Levi-Civita symbol, asserting the expected transformation behavior. This holds leafwise in an asymptotic sense. Consequently, it does not require the STCMC-center of mass $\vec{Z}_{\text{STCMC}}$ to converge. 

The same result holds for asymptotically Schwarzschildean initial data as defined in Section~\ref{sec:AS}. Both results provide an alternative to assuming strong Regge--Teitelboim conditions. They carefully exploit the properties of the STCMC-foliation as well as the Einstein evolution equations and the divergence theorem. They do not rely on a super-potential representation (and such a representation is not known to exist for the STCMC-center of mass).

Thus, the STCMC-center of mass exhibits the expected special-relativistic boost behaviour also under different decay
assumptions than the Regge--Teitelboim conditions, such as $p+q>3$ or asymptotic Schwarzschildeanness. On the other hand, Senthil Velu~\cite{Sara} proves that the B\'ORT-/CMC-center of mass has a \emph{boost defect} in these contexts, that is, a term coming out of the infinitesimal boost law in addition to the Regge--Teitelboim angular momentum. This boost defect depends on the choice of spacetime lapse (or choice of spacetime coordinates), upon boosting. This seems related to the fact that the STCMC-center of mass is defined in a spacetime-covariant fashion and hence does not depend on the spacetime lapse.

\begin{acknowledgement}
 The newer results described in this article were obtained by members of the project \emph{Geometrically defined asymptotic coordinates in general relativity} in the DFG priority program SPP 2026 \emph{Geometry at Infinity}: besides the authors, those are Rodrigo Avalos (Potsdam/T\"ubingen, now Rostock University), Melanie Graf (T\"ubingen/Potsdam, now Hamburg University),  Alejandro Pe\~{n}uela (Potsdam), Annachiara Piubello (Potsdam, now T\"ubingen), Anna Sancassani (T\"ubingen), Saradha Senthil Velu (T\"ubingen), Olivia
  Vi\v{c}\'{a}nek Mart\'{i}nez (T\"ubingen), Markus Wolff (T\"ubingen, now University of Vienna). The listed affiliations listed refer to the timeframe of the project. We would like to cordially thank all project members as well as Domenico Giulini, Gerhard Huisken, and David Maxwell for interesting scientific discussions.
\end{acknowledgement}
\ethics{Competing Interests}{
This research was funded by the Deutsche Forschungsgemeinschaft (DFG, German Research Foundation) -- 441897040.}

\eject

\bibliographystyle{alpha}
 \bibliography{SPP}
\end{document}